\documentclass[11pt]{article}           
\usepackage{a4wide}                     
\usepackage{latexsym}                   
\usepackage{graphicx}
\usepackage{eufrak}

\newfont{\bbold}{msbm10 scaled \magstep1}
\newtheorem{Theorem}{Theorem}
\newtheorem{Proposition}[Theorem]{Proposition}

\newcounter{exampleno}
\newenvironment{Example}{
	\refstepcounter{exampleno}%
 \medbreak\bigskip\noindent%
 {\sc Example~\theexampleno.} }%
 {\nobreak{\hfill$\Box$}\medbreak}
\newcommand{\stir}[2]
{{{#1} \brace {#2}}}
\newcommand{\cfrac}[2]{\frac{\displaystyle#1}{\displaystyle#2}}

\def\Proof{{\noindent \bf Proof. }}

\newcommand{\Ref}[1]{(\ref{#1})}
\def\qed{$\hfill{\vrule height 3pt width 5pt depth 2pt}$ \bigskip}

\def\tilde{\widetilde}

\newcommand{\gf}{generating function}
\newcommand{\gfs}{generating functions}
\newcommand{\beq}{\begin{equation}}
\newcommand{\eeq}{\end{equation}}
\newcommand{\al}{\alpha}
\newcommand{\be}{\beta}
\newcommand{\ga}{\gamma}

\newcommand{\ns}{\mbox{\bbold N}}
\newcommand{\zs}{\mbox{\bbold Z}}

\title{\bf Generating functions\\ for generating trees}

\author{
Cyril Banderier \ \ \ \ \ \ \ \and Mireille Bousquet-M\'elou\thanks
{Corresponding author.}
\ \ \ \ \ \ \  \and Alain Denise  
 \and Philippe Flajolet \and
 Dani\`ele Gardy \and Dominique Gouyou-Beauchamps 
}
\date{}

\begin{document}
\maketitle
\noindent
{\em 
This article corresponds, up to minor typo corrections,
to the article submitted to Discrete Mathematics (Elsevier) in Nov. 1999, 
and published in its vol. 246(1-3), March 2002, pp. 29-55.
This work supplants a 
preliminary draft ``On Generating Functions of Generating Trees'' 
dated Nov. 1998 which appeared in the Proceedings of the FPSAC'99 Barcelona Conference.}

\newcommand{\Sn}{{\EuFrak S}}
\def\abstractfr{\if@twocolumn
\section*{R\'esum\'e}
\else \small 
\begin{center}
{\bf R\'esum\'e\vspace{-.5em}\vspace{0pt}} 
\end{center}
\quotation 
\fi}
\def\endabstractfr{\if@twocolumn\else\endquotation\fi}

\def\abstractangl{\if@twocolumn
\section*{Abstract}
\else \small 
\begin{center}
{\bf Abstract\vspace{-.5em}\vspace{0pt}} 
\end{center}
\quotation 
\fi}
\def\endabstractangl{\if@twocolumn\else\endquotation\fi}

\begin{abstractangl}
Certain families of 
combinatorial objects admit recursive descriptions in terms of
generating trees: each node of the tree corresponds to an object,
and the branch leading to the node encodes
the choices made in the construction of the object. Generating trees
lead to a fast  
computation of enumeration sequences (sometimes, to explicit
formulae as well) and provide efficient random generation algorithms.
We investigate the links between the structural
properties of the rewriting rules defining such trees and the rationality,
algebraicity, or transcendence of the corresponding generating function.
\end{abstractangl}

\section{Introduction}

Only the simplest combinatorial structures --- like binary strings,
permutations, or  pure involutions ({\em i.e.}, 
involutions with no fixed point) ---  admit product
decompositions. In that case, the set $\Omega_n$ of objects of size~$n$
is isomorphic to a product set:
$\Omega_n\cong[1,e_1]\times[1,e_2]\times\cdots\times[1,e_n]$.
Two properties result from such a strong decomposability property:
$(i)$~enumeration is easy, since the cardinality of $\Omega_n$ is
$e_1e_2\cdots e_n$;
$(ii)$~random generation is efficient since it reduces to a sequence
of random independent draws from intervals.
A simple infinite tree, called a {\em uniform generating
tree\/} is determined by the $e_i$: the root has
degree $e_1$, 
each of its $e_1$ descendants has degree~$e_2$, and so on. 
This tree describes the sequence of all possible choices and
the objects of size~$n$
are then in natural correspondence with the branches of length~$n$,
or equivalently with the nodes of generation~$n$ in the tree. The 
generating tree
is thus fully described  by its root degree 
($e_1$) and by rewriting rules, here of the special form, 
\[
(e_i)\leadsto (e_{i+1})\,(e_{i+1})\,\cdots (e_{i+1}) \equiv (e_{i+1})^{e_i},
\]
where the power notation is used to express repetitions. For
instance binary strings, 
permutations, and pure involutions are determined by
\[
\begin{array}{ll}
{\cal S}\,: & [(2),~(2)\leadsto (2)\,(2)] \\
{\cal P}\,: &[(1),~\{(k)\leadsto (k+1)^{k}\}_{k\ge1}]\\
{\cal I}\,:&[(1),~\{(2k-1)\leadsto (2k+1)^{2k-1}\}_{k\ge1}].
\end{array}
\]

A powerful generalization of this idea consists in
considering unconstrained {\em generating trees\/} where any set of rules
\begin{equation}\label{i1}
\Sigma=[(s_0),~\left\{(k)\leadsto
(e_{1,k})\,(e_{2,k})\,\cdots\,(e_{k,k})\right\}]
\end{equation}
is allowed. Here, the {\em axiom\/} $(s_0)$ specifies the degree of
the root,
while the {\em productions\/} $e_{i,k}$ list the degrees of the $k$ descendants
of a node labeled $k$. 
Following Barcucci, Del~Lungo, Pergola and
Pinzani, we call $\Sigma$ an
{\em ECO-system\/} (ECO stands for ``Enumerating
Combinatorial Objects'').
Obviously, much more leeway is available and there is hope
to describe a much wider class of structures than those corresponding
to product forms and uniform generating trees.

The idea of generating trees 
has surfaced
occasionally in the literature. West introduced it in the context of
enumeration of permutations with forbidden subsequences~\cite{West95,West96};
this idea has been further exploited  in closely related
problems~\cite{BaDePePiFPSAC98,BaDePePi,DuGiGu98,DuGiWe96}. A major
contribution in this area is due to Barcucci, Del~Lungo, Pergola, and
Pinzani~\cite{ECOArticle,BaDePePi_2} who 
showed that a fairly large number of
classical combinatorial structures can be described 
by generating trees. 

A form equivalent to generating trees is well worth noting at this
stage. Consider the  {\em walks\/} on the integer half-line that start
at~point $(s_0)$ and such that the only allowable transitions are 
those specified by~$\Sigma$ {(the steps corresponding to
transitions with multiplicities being labeled)}. Then, the walks of
length~$n$  
are in bijective correspondence with the nodes of generation $n$ in the
tree. 
These walks are constrained by the
consistency requirement of trees, namely, that the number of  outgoing
edges from point $k$ on the half-line has to be exactly~$k$.

\begin{figure}[ht]
\begin{center}
\includegraphics{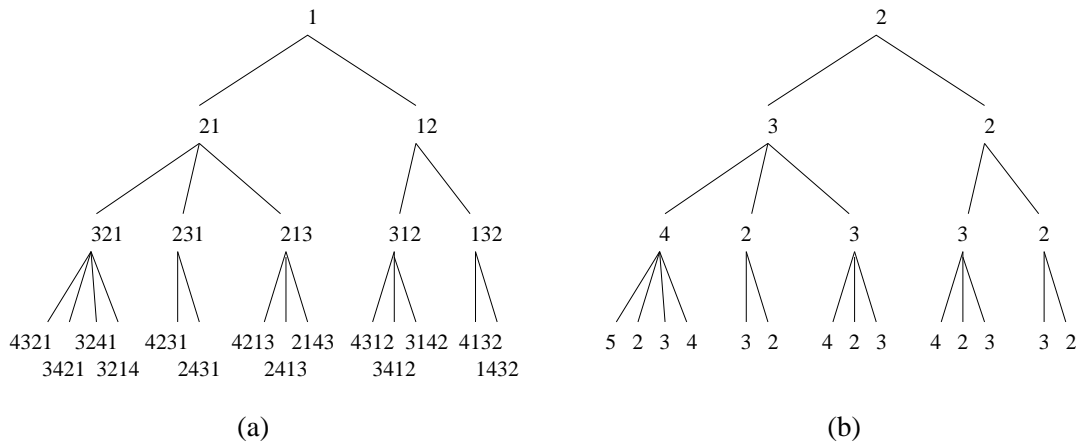}
\end{center}
\caption{The generating tree of $123$-avoiding permutations. 
(a)~Nodes labeled by the permutations. 
(b)~Nodes labeled by the numbers of children.}
\label{pinzani-arbre1}
\smallskip
\noindent\hrule
\end{figure}

\begin{Example} {\em $123$-avoiding permutations}\label{Catalan123Example}\\
The method of ``local expansion'' sometimes gives good results in
the enumeration of permutations avoiding specified patterns.
 Consider for example the set $\Sn _n(123)$ of permutations  of length
$n$ that {\em avoid the pattern\/} $123$:
there exist no integers $i < j< k$ such that $\sigma(i) <
\sigma(j) < \sigma(k)$. For instance, $\sigma=4213$ belongs to
$\Sn _4(123)$ but $\sigma=1324$ does not, as $\sigma(1) < \sigma(3)
< \sigma(4)$.

Observe that if $\tau \in \Sn _{n+1}(123)$, then the permutation $\sigma$
obtained by erasing the entry $n+1$ from $\tau$ belongs to $ \Sn
_{n}(123)$.
Conversely, for every
$\sigma \in \Sn _{n}(123)$, insert the value $n+1$ in each  place that gives an element of $\Sn _{n+1}(123)$
(this is the local expansion). For example, the permutation
$\sigma=213$ gives $4213$, $2413$ and $2143$, by insertion of $4$ in
first, second and third place respectively. The permutation $2134$,
resulting {from} the insertion of $4$ in the last place, does not belong
to $\Sn _4(123)$. This process can be described by a 
tree whose nodes are the permutations avoiding $123$:
the root is $1$, and  the children of any node $\sigma$ are the
permutations derived as above. Figure~\ref{pinzani-arbre1}(a) presents
the first four levels of this tree.  

 Let us now label the nodes by their number of children: we obtain the
 tree of Figure~\ref{pinzani-arbre1}(b). It can be proved that the 
$k$ children of any node labeled $k$ are labeled respectively
$k+1,2,3,\ldots,k$ {(see \cite{West95})}. Thus the  tree we have
 constructed is the 
generating tree obtained from the following rewriting rules: 
$$
	[ (2),\ \{(k)  \leadsto (2)(3)\ldots(k-1)(k)(k+1) \}_{k\ge 2}].
$$
The interpretation of {this system} in terms of paths implies that $123$-avoiding permutations  
are equinumerous with ``walks with returns'' on the half-line,
themselves isomorphic to {\L}ukasiewicz codes of plane trees (see,
{\em e.g.}, \cite[p.~31--35]{stanley-vol2}).
We thus recover a classic result \cite{Knuth1}:
 $123$-avoiding permutations are counted by  Catalan numbers; more precisely,
$|\Sn_n(123)|={2n \choose n}/(n+1)$.
\end{Example}

We shall see below that (certain) generating trees correspond to
enumeration sequences of relatively low computational 
complexity and provide  fast
random generation  algorithms. Hence, there is an obvious interest in
delineating as precisely as 
possible which combinatorial classes admit a generating tree
specification. Generating functions condense structural
information in a simple analytic entity. We can thus wonder
what kind of \gf \ can be obtained through generating trees. More
precisely, we study in this paper  the
connections between the {\em structural\/} properties of the rewriting rules
and the {\em algebraic\/} properties of the corresponding generating function.

We shall prove several conjectures that were presented to
us by Pinzani and his coauthors in March 1998.
 Our main results can be roughly described as follows.
\begin{itemize}
\item[---] {\em Rational systems.} Systems satisfying strong
regularity conditions
lead to rational generating functions (Section~\ref{SectionRational}).
This covers systems that
have a finite number of allowed degrees, as well as 
systems like $(2.a)$, $(2.b)$, $(2.c)$ and $(2.d)$
 below where the labels are constant except for 
a fixed number of labels that depend linearly  and uniformly on $k$.
\item[---] {\em Algebraic systems.} 
Systems of a {\em factorial\/} form,
{\em i.e.}, where a finite modification of the set $\{1,\ldots,k\}$ is
reachable 
from~$k$, lead to algebraic generating functions
(Section~\ref{SectionAlgebraic}). This includes in particular
cases  $(2.f)$ and $(2.g)$.
\item[---] {\em Transcendental systems.} One possible reason for a
system to give a transcendental series is the 
fact that its coefficients grow too fast, so that its radius of convergence is
zero. This is the case for System~$(2.h)$ below.
Transcendental generating functions are also associated with systems
that are too ``irregular''. An example is System~$(2.e)$.
 We shall also discuss the
holonomy of transcendental systems (Section~\ref{SectionTranscendant}).
\end{itemize}

\newpage
\begin{Example} \label{SystemExamples}
{\em A zoo of rewriting systems}\\
Here is a list of
examples  recurring throughout this paper.
\begin{flushright}
\begin{tabular}{ll}
$ [(3), \{(k) \leadsto (3)^{k-3} (k+1) (k+2) (k+9)\}]$
  & $(2.a)$  \\
$ [(3), \{(k) \leadsto (3)^{k-1} (3k+6)\}]$
  & $(2.b)$ \\
$ [(2), \{(k)  \leadsto (2)^{k-2} (2+(k \bmod 2)) (k+1)\}]$
  & $(2.c)$ \\
$ [(2), \{(k)  \leadsto (2)^{k-2} (3-(k \bmod 2)) (k+1)\}]$
  & $(2.d)$ \\
$ [(2), \{(k) \leadsto (2)^{k-2} (3 - [\exists p\!:\!k=2^p]) (k+1)\}]$
 \ \ \ \ \ \ \ \ \ \ \ \ \ \ \ \ \ \ \ \  & $(2.e)$ \\
$ [(2), \{(k) \leadsto (2)(3)\ldots(k-1)(k)(k+1)\}]$
  & $(2.f)$ \\
$ [(1), \{(k) \leadsto (1)(2)\ldots(k-1) (k+1)\}]$
  & $(2.g)$ \\
$ [(2), \{(k) \leadsto (2)(3) (k+2)^{k-2}\}]$ &$(2.h)$ \\
\end{tabular}
\end{flushright}
(In $(2.e)$, we make use of Iverson's brackets: $[P]$ equals $1$ if $P$ is
true, $0$ otherwise.)
\end{Example}

\paragraph{Notations.}
{}From now on, we adopt functional notations for rewriting rules:
systems will be of the form
$$
  [(s_0), \quad \{(k) \leadsto (e_{1}(k))\; (e_{2}(k)) \ldots (e_{k}(k))\}]
$$
where $s_0$ is a constant and each $e_i$ is a function of
$k$. Moreover, we assume that 
all the values appearing in the generating tree are positive: {each
node has at least one descendant.}

In the generating tree, let $f_n$ be the number of nodes at level $n$
and $s_n$ the sum of the labels of these nodes. By convention, the
root is at level $0$, so that $f_0=1$. In terms of walks, $f_n$ is the
number of walks of length $n$.
The generating function associated {with} the system is
$$
	F(z) = \sum_{n \geq 0} f_n z^n.
$$
Remark that $s_n=f_{n+1}$, and that the sequence $(f_n)_n$ is
nondecreasing.

Now let $f_{n,k}$ be the number of nodes at level $n$ having
label $k$ (or the number of walks of length $n$ ending at position
$k$). The following generating functions will be also of interest:
$$
F(z,u) = \sum_{n,k \geq 0} f_{n,k} z^n u^k
\qquad \hbox{and} \ \ 
F_k(z) = \sum_{n \geq 0} f_{n,k} z^n.
$$
We have $F(z)=F(z,1)=\sum_{k\ge 1}F_k(z)$.
Furthermore, the $F_k$'s satisfy the relation
\begin{equation} \label{equationF_k}
F_k(z) = [k=s_0] + z \sum_{j \geq 1} \pi_{j,k} F_j(z),
\end{equation}
where $\pi_{j,k} = |\{i \le j: e_i(j)=k \}|$ denotes the number of one-step
transitions from $j$ to $k$. This is equivalent to the following
recurrence for the numbers $f_{n,k}$,
\begin{equation}\label{rec1}
{f_{0,k}=[k=s_0]}  \qquad 
	\hbox{and} \ \ f_{n+1,k}=\sum_{j \geq 1} \pi_{j,k} f_{n,j},
\end{equation}
that results from tracing all the paths that lead to $k$ in $n+1$
steps.  

\paragraph{Counting and random generation.}
The recurrence~(\ref{rec1})
gives rise to an algorithm that computes  the successive rows of the
matrix  $(f_{n,k})$ by ``forward propagation'': to compute the
$(n+1)$th row, propagate the contribution $f_{n,j}$ to
$f_{n+1,e_i(j)}$ for all pairs $(i,j)$ such that $i\le j$.
Assume the system is  {\em linearly bounded\/}: this means that the
labels of the nodes that can be  reached in $m$ steps are bounded by
a linear function of~$m$. (All the systems given in
Example~\ref{SystemExamples}, except for $(2.b)$, are  
linearly bounded; more generally, systems where forward jumps are bounded by a
 constant are linearly bounded.)
Clearly, the forward propagation algorithm 
 provides a counting algorithm of arithmetic  complexity that is at most
cubic. 

For a linearly bounded system, uniform random generation can also be achieved
in polynomial 
time, as shown in \cite{BaDePe99}. We present here the general
principle.

Let $g_{n,k}$ be the number of walks of length~$n$ that start from
label~$k$.
These numbers
are determined by  the  recurrence
$g_{n,k}=\sum_i g_{n-1,e_i(k)},$
that traces all the possible continuations of a path given its initial
step. Obviously, $f_n=g_{n,s_0}$, with $s_0$ the axiom of the system.
As above, the $g_{n,k}$ 
can be determined in time $O(n^3)$ and~$O(n^2)$ storage.
Random generation is then achieved as follows:
In order to generate a walk of length $n$ starting from state~$k$,
pick up a transition $i$ with probability
$g_{n-1,e_i(k)}/g_{n,k}$,
and generate recursively a walk of length~$n-1$ starting from state
$e_i(k)$. 
The cost of a 
single random generation is then  $O(n^2)$ if a sequential search is
used over the $O(n)$ possibilities of each of the $n$ random drawings;
the time complexity goes down to $O(n\log n)$ if binary search is
used, but at the expense of an increase in storage complexity of
$O(n^3)$ (arising from $O(n^2)$
arrays of size $O(n)$ that binary search requires).

\section{Rational systems}
\label{SectionRational}

We give in this section three main criteria (and a variation on one of them)  
implying that the generating
function of a given ECO-system is rational.

Our first and simplest criterion applies to systems in which
 the functions $e_i$ are  uniformly bounded.
\begin{Proposition}\label{obvious}
If finitely many labels
appear in the tree, then $F(z)$ is rational.
\end{Proposition}

\Proof
Only a finite number of $F_k$'s are nonzero, and they are related by
linear equations like Equation~(\ref{equationF_k}) above.
\qed

\begin{Example} \label{FibonacciExample}
{\em The Fibonacci numbers} \\
The system $[(1),\, \{(k) \leadsto (k)^{k-1}
((k \bmod 2)+1)\}]$ can be also written as $[(1),\, \{(1)
\leadsto (2),\, (2) \leadsto (1) (2)\}]$. Hence the only labels that
occur in the tree are $1$ and $2$. Eq.~\Ref{equationF_k} gives $F_1(z)=1+zF_2(z)$
and $F_2(z)=z(F_1(z)+F_2(z))$. Finally,
$$F(z)=\frac{1}{1-z-z^2}=\sum_{n\ge 0} f_n z^n= 1+z+2z^2+3z^3+5z^4+\cdots,$$ 
the well-known Fibonacci generating function.
\end{Example}

None of the systems of Example~\ref{SystemExamples} satisfy  the
assumptions of Proposition~\ref{obvious}.  However, the following criterion can
be applied to systems $(2.a)$ and~$(2.b)$. 

\begin{Proposition} \label{Prop1Rule}
Let $\sigma(k)=e_1(k)+e_2(k)+\cdots +e_k(k)$. 
If $\sigma$ is an affine function of $k$, say $\sigma(k)=\alpha k+\beta$, then
the series $F(z)$ is rational. More precisely:
$$
F(z) \ =\  {1 + (s_0-\alpha)z \over 1-\alpha z-\beta z^2}.
$$
\end{Proposition}
\Proof
Let $n \geq 0$ and let $k_1, k_2, \ldots k_{f_n}$ denote the labels of
the $f_n$ nodes at level $n$. Then
\begin{eqnarray*}
f_{n+2} \ =\  s_{n+1} &=& (\alpha k_1+\beta) + (\alpha k_2+\beta) 
					+ \cdots + (\alpha k_{f_n}+\beta) \\
	          &=& \alpha s_n + \beta f_n \ =\  \alpha f_{n+1} + \beta f_n.
\end{eqnarray*}
We know that $f_0=1$ and $f_1=s_0$. The result follows. 
\qed

\begin{Example} \label{FiboBisectionExample}
{\em Bisection of Fibonacci sequence}  \\
The system $[(2),\,\{(k)
\leadsto (2)^{k-1} (k+1)\}]$ gives $F(z)={1-z \over
1-3z+z^2}=1+2z+5z^2+\cdots$,  
the generating function for Fibonacci numbers of even index.
(Changing the axiom to $(s_0)=(3)$ leads to the other half of the Fibonacci
sequence.) 
Some other systems, like
$$\begin{array}{l}
{[(2),\,\{(k) \leadsto (1)^{k-1} (2k)\}]}, \\
{[(2),\,\{(k) \leadsto (2)^{k-2} (3-(k \bmod 2)) (k + (k \bmod 2))\}]} ,\\
{[(2),\,\{(k) \leadsto (2)^{k-2} (3-[k \mbox{ is prime}]) (k + [k \mbox{ is prime}])\}],}
\end{array}$$
lead to the same function $F(z)$ since $\sigma(k)=3k-1$  
and $s_0=2$. However, the generating trees are
different, as are the bivariate functions $F(z,u)$.
\end{Example}

\begin{Example} \label{ParodoxExample}
{\em Prime numbers and rational \gfs}\\
Amazingly, it is  possible to construct a
generating tree whose set of labels is the set of prime
numbers but that has a rational generating function $F(z)$.
This is a bit unexpected, as prime numbers are usually thought ``too
irregular'' to be associated with rational \gf s.
For $n\ge 1$, let $p_n$ denote the $n$th prime; hence $(p_1, p_2, p_3,
\ldots) = (2,3,5,\ldots)$. Assume for the moment that the Goldbach
conjecture is true: every even number larger than $3$ is the sum of
two primes. Remember that, according to  Bertrand's postulate, 
$p_{n+1}< 2p_n$ for all $n$ (see, {\em e.g.}, \cite[p.~140]{Ribenboim}).

For $n \ge 1$, the number $2p_n-p_{n+1}+3$ is an even number larger
than $3$. Let $q_n$ and $r_n$ be two primes such that
$2p_n-p_{n+1}+3=q_n+r_n$. In particular, $q_1=r_1=2$. Consider the
system
$$[(2), \{ (p_n)\leadsto (p_{n+1}) (q_n) (r_n) (2)^{p_n-3}\}].$$
It satisfies the criterion of Proposition~\ref{Prop1Rule}, with
$\sigma(k)=4k-3$. Hence, the \gf \ of the associated generating tree
is
$$F(z)=\frac{1-2z}{1-4z+3z^2}=\frac{1}{2} \left[
\frac{1}{1-z}+\frac{1}{1-3z}\right].$$
Consequently, the number of nodes at level $n$ is simply {$f_n=(1+3^n)/2$.}
This can be checked on the first few levels of the tree drawn in
Figure~\ref{goldbach}. 

\begin{figure}[ht]
\begin{center}
\includegraphics{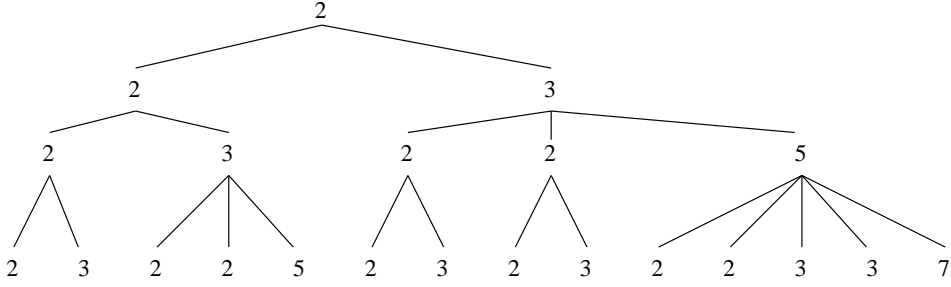}
\end{center}
\caption{A generating tree with prime labels and rational \gf .}
\label{goldbach}
\smallskip
\end{figure}

Now, one can object that the Goldbach conjecture is not proved;
however, it is known that every even number is the sum of at most six
primes~\cite{Ramare}, and a similar example can be constructed using
this result. 

\end{Example}

\medskip
Proposition~\ref{Prop1Rule} can be adapted to apply to systems that
``almost'' satisfy the criterion of Proposition~\ref{Prop1Rule}, like
   System~$(2.c)$ {or~$(2.d)$.}
Let us consider a system of the form
\begin{eqnarray*}
(s_0), & (k) \leadsto e_1^{[0]}(k), \ldots , e_k^{[0]}(k) \ \ \hbox{ if }
k \hbox{ is  even},\\
& (k) \leadsto e_1^{[1]}(k), \ldots , e_k^{[1]}(k) \ \ \hbox{ if }
k \hbox{ is  odd}.
\end{eqnarray*}
Assume, moreover, that:

$(i)$~the corresponding functions $\sigma_0$ and $\sigma_1$ are
affine and have the same leading coefficient $\alpha$, say $\sigma_0(k)=\alpha
k+\beta_0$ and $\sigma_1(k)=\alpha k+\beta_1$;

$(ii)$~exactly $m$
odd labels occur in the right-hand side of each rule, for some $m\ge 0$.
\begin{Proposition} \label{Prop2Rules}
If a system satisfies properties {\em ({\it i})} and {\em ({\it ii})}
 above, then
$$
F(z) = {1 + (s_0-\alpha)z + (s_1-\alpha s_0-\beta_0)z^2
       \over
	1 - \alpha z - \beta_0z^2 - m(\beta_1-\beta_0)z^3}.
$$
Of course, if $\beta_0=\beta_1$, we recover the \gf \ of
Proposition~{\em \ref{Prop1Rule}. } 
\end{Proposition}
\Proof The proof is similar to that of
       Proposition~\ref{Prop1Rule}. The only new ingredient is the
       fact that, for 
       $n\ge 1$, the number of nodes of odd label at level $n$ is
       $mf_{n-1}$. 
\qed

 System~$(2.c)$ satisfies properties
({\it i}) and ({\it ii}) above with 
$\alpha=3$, $\beta_0=-1$, $\beta_1=0$, $m=1$, $s_0=2$ and $s_1=5$. 
Consequently, its generating
function is $F(z)= { 1-z \over 1-3z+z^2-z^3}$.
System~$(2.d)$, although  very close to $(2.c)$,
does not satisfy property ({\it ii}) above, so that
Proposition~\ref{Prop2Rules} does not apply.  
{However, another minor variation on the argument of
Proposition~\ref{Prop1Rule}, based on the fact that the number $o_n$ of
odd labels at level $n$ satisfies $o_n=2(f_{n-1}-o_{n-1})$, proves the
rationality of $F(z)$.}

{Alternatively, rationality follows from the last 
criterion of this section, which is of a different nature.}
We  consider systems $[(s_0),\,
\{(k) \leadsto (e_1(k)) (e_2(k)) \ldots (e_{k}(k))\}]$ that can be
written as
\begin{equation}
[(s_0),\,
\{(k) \leadsto (c_1(k)) (c_2(k)) \ldots (c_{k-m}(k))
			(k+a_1)
				(k+a_2)\ldots
					(k+a_m)\}]
\label{lineaire1}
\end{equation}
where $1\le a_1\le a_2\le \cdots\le a_m$ and
the functions $c_i$ are uniformly  bounded.
Let $C =\max_{i,k} \{s_0,  c_i(k)\}$.

\begin{Proposition}
Consider the system {\em (\ref{lineaire1})},  and  let 
$\pi_{j,k} = |\{i\le j : e_i(j)=k \}|$. If all  the series
$$
\sum_{j \ge 1} \pi_{j,k} \ t^j
$$
for $k \le C$ are rational, then so is the series $F(z)$.
\end{Proposition}

\Proof
We form an infinite system of equations defining the series $F_k(z)$ by writing
Eq.~(\ref{equationF_k}) for all $k\ge 1$. In particular, for 
$k> C$, we obtain
$$F_k(z)=z\sum_{\ell=1}^m  F_{k-a_\ell}(z),$$
with $F_j(z)=0$ if $j\le 0$.
This part of the system  is easy to solve in terms of  $F_1, \ldots
,F_C$. Indeed, for $k\in \zs$: 
\begin{equation}
\label{recFP}
	F_k(z) = \sum_{i=1}^{C}P_{i,k}(z) F_i(z)
\end{equation}
where the $P_{i,k}$ are polynomials in $z$
defined by the following recurrence: for all $i\le C$,
\begin{equation}
P_{i,k}(z) = \left\{
\begin{array}{ll}
  0   & \hbox{\rm if } k \le 0,  \\
{[k=i]} & \hbox{\rm if } 0<k \le C ,\\	
 \displaystyle
z \sum_{\ell=1}^m{ P_{i,k-a_\ell}(z)} &
			\mbox{\rm if } k >C. \\
\end{array}\right.
\label{Pik}
\end{equation}
\medskip
Using (\ref{recFP}), we find
$$F(z)=\sum_{k\ge 1}F_k(z)= \sum_{i=1}^C\left[ F_i(z)  \sum_{k \ge 1}
P_{i,k}(z) \right].$$
According to (\ref{Pik}), for all $i\le C$, the series 
$\sum_{k \ge 1} P_{i,k}(z)t^k$ is a rational
function of $z$ and $t$, of denominator
$1-z\sum_\ell t^{a_\ell}$. At $t=1$, it is rational in $z$.
Hence,  to prove the rationality of $F(z)$, it suffices to prove the
rationality  of the $F_i(z)$, for $i\le C$.

Let us go back to  the $C$ first equations of our system; {using
(\ref{recFP})}, we find, for $k\le C$:
$$
	F_k(z) = [k=s_0] + z \sum_{i=1}^C\left[ F_i(z)\sum_{j \ge 1}P_{i,j}(z)\pi_{j,k}\right].
$$
Again,  $\sum_{j \ge 1}P_{i,j}(z)\pi_{j,k}t^j$ is a
rational function of $z$ and $t$ (the Hadamard product of two rational series
is rational). Thus the series $F_k(z)$, for $k\le C$, satisfy a linear
system with rational coefficients: they are rational themselves, as
well as $F(z)$.
\qed

Examples~$(2.a)$, $(2.c)$, $(2.d)$ and $(2.e)$
have the form (\ref{lineaire1}). The above proposition implies that
the first three have a rational generating function.  System~$(2.e)$
will be discussed in Section~\ref{SectionTranscendant}, and proved to
have a transcendental \gf . 

\section{Factorial walks and algebraic systems}
\label{SectionAlgebraic}

In this section, we consider systems that are of a {\em factorial
form\/}. By this, we  mean informally that the set of
successors of $(k)$ is a finite modification of the integer interval
$\{1,2,\ldots,k\}$. As was 
detailed in the introduction, ECO-systems can be rephrased
in terms of walks over the integer half-line. We 
thus consider the 
problem of enumerating walks
over the integer half-line such that the set of allowed moves from
point~$k$ is a finite modification  of the
integer interval 
$[0,k]$. We shall mostly study modifications around the point $k$
(although some examples where the interval is modified around $0$ as
well are given at the end of the section).
Precisely, a {\em factorial walk\/} is
defined by a finite (multi)set $A \subset \zs$ and a finite set $B \subset \ns^+$, where $\ns ^+=\{1,2,3, \dots\}$,
 specifying respectively the {\em
allowed supplementary jumps\/} (possibly labeled)
and the {\em forbidden backward jumps\/}. In
other words, the possible moves from $k$ are given by the rule:
\begin{equation}\label{alg1}
(k)\leadsto [0,k-1] \setminus (k-B) \ \cup \  (k+A).
\end{equation}
Observe that these walk models are not necessarily ECO-systems, first
because we allow labels to be zero -- but a simple translation can
take us back to a model with positive labels -- and second because we
do not require $(k)$ to have exactly $k$ successors.

We say that an ECO-system is factorial if a shift of indices
transforms it into a factorial walk. Hence the rules of a factorial
ECO-system are of the form 
\[
(k+r)\leadsto [r, k+r-1] \setminus (k+r-B) \ \cup \ (k+r+A),
\]
that is, 
\beq \label{eco-factorial}
(k)\leadsto 
[r,k-1] \setminus(k-B) \ \cup \ (k+A) \ \ \ \hbox{ for } k \ge r \ge 1.
\eeq
The \gf \ $F(z)$ for
such an ECO-system, taken with axiom $(s_0)$, equals the \gf \ for the
walk model \Ref{alg1}, taken with axiom $(s_0-r)$. However, remember
that 
the rewriting rules defining a generating tree have to obey the
additional condition that a node labeled $k$ has exactly $k$
successors. Taking $k=r$ in \Ref{eco-factorial}, this implies that $r=|A|$. Taking
$k>r+\max B$, this implies that $r+|B|=|A|$, so that finally $B =
\emptyset$. Hence, strictly speaking, either 
one has a ``fake'' factorial ECO-system (that is some of its initial
rules are not of the factorial type), 
either one has a ``real'' factorial ECO-system and then it is given by rules of the form
\[
(k)\leadsto 
[r,k-1]  \ \cup \ (k+A) \ \ \ \hbox{ for } k \ge r \ge 1,
\]
where $A$ is a multiset of integers of cardinality $r$.
For instance, Systems~$(2.f)$ and $(2.g)$ are factorial.
We shall prove that
all factorial  walks have an algebraic generating function.
The result naturally applies to factorial ECO-systems.

We consider again the generating function
$F(z,u)=\sum_{n,k\geq 0}f_{n,k} z^n u^k$, where 
$f_{n,k}$ is the number of walks of length $n$  ending at  point
$k$. We also denote by $F_k(z)$  the coefficient of $u^k$ in this series,
and by $f_n(u)$  the coefficient of $z^n$.
The first ingredient of the proof is {a} linear operator $M$, acting on
formal power series in $u$, that encodes the possible moves. More
precisely, for all $n\ge 0$, {we will have:}
\[
M[f_n](u)=f_{n+1}(u).
\]
The operator $M$ is constructed step by step as follows.
\begin{itemize}
\item[---]
The set of moves from $k$ to all the positions $0, 1,
\dots, k-1$ is described by the operator $L_0$ that maps $u^k$ to 
$u^0+u^1+\cdots+u^{k-1} =(1-u^k)/(1-u)$. 
As $L_0$ is a linear operator, we have, for any series $g(u)$:
\[
L_0[g](u)=\frac{g(1)-g(u)}{1-u}.
\]
\item[---]
The fact that transitions near $k$ are modified, with those of type
$k+\al$ (with $\al \in A$)
allowed and those of type $k-\be $ (with $\be \in B $)
forbidden, is expressed by a Laurent polynomial
$$
P(u)=\sum_{k=-b}^a p_k u^k=A(u)-B(u)\quad\hbox{with}\quad
A(u)=\sum_{\al \in A} u^{\al} \ \ \ \hbox{and }\ \ \ 
B(u)=\sum_{\be \in B} u^{-\be}.
$$
The degree  of $P$ is $a:= \max A$, the largest forward jump
and $b:=\max(0,-B,-A)$ is largest forbidden backward jump or the largest
  supplementary backward jumps
(we take $b=0$ if the set $B$ is empty). 

The operator 
$$L[g](u):=L_0[g](u) + P(u)g(u)$$
describes the extension of a walk by one step.
\item[---] Finally, the operator $M$ is given by
\[
M[g](u)=L[g](u)-\{u^{<0}\} L[g](u),
\]
where $\{u^{<0}\}h(u)$ is the sum of all the monomials in
$h(u)$ having a
negative exponent. Hence $M$ is nothing but $L$ stripped of the negative
exponent monomials, which correspond to walks ending on the nonpositive
half-line. Observe that, for any series $g(u)$, the only part of
$L[g](u)$ that is likely to contain monomials with negative exponents
is $P(u)g(u)$. Consequently, 
\[
M[g](u)=L[g](u)-\{u^{<0}\}[P(u)g(u)]
\]
and if $g(u)=\sum_k g_k u^k$, then
\begin{equation}
\{u^{<0}\}[P(u)g(u)]
=\sum_{i=1}^{b}  \sum_{k=0}^{i-1}g_k p_{-i} \,  u^{k-i}
=\sum_{k=0}^{b-1} g_k r_k(u). \label{<0-simple}
\end{equation}
\end{itemize}
Assume for simplicity that the
initial point of the walk is~$0$;
other cases follow  the same argument.
The linear relation $f_{n+1}(u)=M[f_n](u)$, together with
$f_0(u)=1$, yields 
\begin{eqnarray*}
F(z,u) &=&1+zM[F](z,u)\\
&=& 1+ z \left( \frac{F(z,1)-F(z,u)}{1-u}+ P(u) F(z,u) +
\{u^{<0}\}  [P(u)F(z,u)] \right).
\end{eqnarray*}
Thanks to \Ref{<0-simple}, we can write
$$\{u^{<0}\}  [P(u)F(z,u)]=
\sum_{k=0}^{b-1} r_k(u)  F_k(z),$$
 where $r_k(u)$ is a Laurent polynomials (defined by Equation~\ref{<0-simple}) 
whose exponents belong to $[k-b, -1]$.
Thus, $F(z,u)$ satisfies the following  functional  equation: 
\begin{equation} \label{e1} F(z,u) \left(1+\frac{z}{1-u}-z P(u)\right)= 
1 +\frac{zF(z,1) }{1-u} 
+z \sum_{k=0}^{b-1} r_k(u)  F_k(z).  \end{equation}

\medskip
Let us take an example. The moves
\[
(k)\leadsto (0)(1)\cdots(k-5) (k-3)(k-1)(k)(k+7)(k+9),
\]
lead to  $A(u)=u^0+u^7+u^9$ and $B(u)=u^{-4}+u^{-2}$. Moreover,
 $$\{u^{<0}\}  [B(u)F(z,u)]=(u^{-2}+u^{-4})F_0(z)
+(u^{-1}+u^{-3})F_1(z)+u^{-2}F_2(z)
+u^{-1} F_3(z),$$
so that the functional equation defining $F(z,u)$ is 
$$ F(z,u) \left(1+\frac{z}{1-u}-z (1+u^7+u^9-u^{-4}-u^{-2})\right)= 
\hskip 50mm {} $$
$$1 +\frac{zF(z,1) }{1-u}  +
z(u^{-2}+u^{-4})F_0(z)
+z(u^{-1}+u^{-3})F_1(z)+zu^{-2}F_2(z)
+zu^{-1} F_3(z).$$

\bigskip

The second ingredient of the proof,  sometimes called the
{\em kernel method\/}, seems to belong to the
``mathematical folklore'' since the 1970's. It has been used in
various combinatorial problems~\cite{CoRi,Knuth1,Odlyzko}
and in probabilities~\cite{FaIa}. See also~\cite{Bo98,BoPe98,Pe98} for
more recent and systematic applications.
This method consists in cancelling
the left-hand side  of
the fundamental functional equation~(\ref{e1}) 
by coupling $z$ and $u$, so that the coefficient of 
the (unknown) quantity $F(z,u)$ is zero. This constraint defines $u$ as one of
the branches of an algebraic
function of $z$. Each branch that can be
substituted analytically into the functional equation yields a linear
relation between the unknown series $F(z,1)$ and 
$F_k(z)$, $0\le k <b$. If enough branches can be
substituted analytically, 
we obtain a system of linear equations, whose solution gives
$F(z,1)$ and the $F_k(z)$  as algebraic
functions. {}From there, an expression for $F(z,u)$ also results in the
form of a bivariate algebraic function.

Let us multiply Eq.~(\ref{e1}) by $u^b(1-u)$ to obtain an equation
with polynomial coefficients {(remind that we take $b=0$ if the set $B$ of
forbidden backward steps is empty)}. The new equation reads
$K(z,u)F(z,u)=R(z,u)$, where $K(z,u)$ is the {\em kernel\/} of the equation: 
\begin{eqnarray}
K(z,u)&=&
u^b(1-u)\left(1+\frac{z}{1-u}-z P(u)\right),\nonumber\\
&=& u^b(1-u)+zu^b-z(1-u)\sum_{\al \in A} u^{\al+b}+z(1-u)\sum_{\be \in
B}  u^{b-\be}.\label{kern}
\end{eqnarray}
This polynomial has degree  $a+b+1$ in $u$, and hence, admits $a+b+1$
solutions, which are algebraic functions of $z$.
The classical theory of algebraic functions
and the Newton polygon construction enable us to expand the
solutions near any point as Puiseux series
(that is, series involving fractional exponents; see~\cite{Dieudonne}).
The $a+b+1$ solutions, expanded around 0,  can be classified as follows:
\begin{itemize} 
\item[---] the ``unit'' branch, denoted by $u_0$,  {is a power
series in $z$ with constant term $1$;}
\item[---] $b$ ``small'' branches, denoted by  $u_1, \dots, u_{b}$, 
{are power series in $z^{1/b}$ whose first nonzero term is $\zeta z^{1/b}$,
with $\zeta ^b+1=0$;}
\item[---] $a$ ``large'' branches, denoted by $v_1, \dots, v_a$, 
{are Laurent series in $z^{1/a}$ whose first nonzero term is $\zeta
z^{-1/a}$, with $\zeta ^a+1=0$.} 
\end{itemize}
{In particular}, all the roots are distinct. {(It is not difficult
to check ``by hand'' the existence of these solutions: for instance,
plugging $z=t^b$ and $u=tw(t)$ in $K(z,u)=0$ confirms the existence of
the $b$ small branches.)}
Note that there are exactly $b+1$ finite branches: the unit
branch $u_0$ and the $b$ small branches $u_1,\ldots,u_b$.
{As $F(z,u)$ is a series in $z$ with \em{polynomial}} {coefficients
in $u$, these $b+1$ series $u_i$, having no negative exponents, can be
substituted for $u$ in $F(z,u)$. More specifically, let us replace $u$
by $u_i$ in~(\ref{e1}): the right-hand side of the equation vanishes,
giving a linear equation relating the $b+1$ unknown series $F(z,1)$ and
$ F_k(z)$, $0\le k <b$. Hence the $b+1$ finite branches give a set of
$b+1$  linear equations relating the $b+1$ unknown series.  One could solve directly this
system, but the following argument is more elegant.}

The right-hand side of~(\ref{e1}), once multiplied by $u^b(1-u)$, is
$$R(z,u)= u^b (1-u) \left(1+\frac{z}{1-u} F(z,1)
+z \sum_{k=0}^{b-1} r_k(u)  F_k(z)\right).$$
By construction, it is a  {\em polynomial\/} in~$u$ of degree $b+1$ and
leading coefficient $-1$. {Hence, it admits $b+1$ roots, which
depend on $z$.} Replacing $u$ by  the series $u_0, u_1, \dots, u_b$ in
Eq.~\Ref{e1} shows that these series are exactly the $b+1$
roots of $R$, so that
$$R(z,u)=-\prod_{i=0}^b (u-u_i).$$
Let $p_a:=[u^a]P(u)$ 
be the multiplicity of the largest forward jump.
Then the  coefficient of $u^{a+b+1}$ in $K(z,u)$ is $p_az$,
and we can write
$$K(z,u)=  p_a z  \prod_{i=0}^{b} (u-u_i)\prod_{i=1}^{a} (u-v_i).$$
Finally, as $K(z,u)F(z,u)=R(z,u)$, we obtain
\begin{equation}\label{biv}
F(z,u) =\frac{-\prod_{i=0}^{b} (u-u_i)}
{u^b(1-u)+zu^b-zu^b(1-u)P(u)}
=-\frac{1}{ p_a z \prod_{i=1}^{a} (u-v_i)}
. 
\end{equation}
We have thus proved the following result.
\begin{Proposition}\label{thm-algebraic}
The \gf \ $F(z,u)$ for factorial walks defined by {\em \Ref{alg1}} and
starting from $0$  is algebraic; it is given by {\em \Ref{biv}}, where 
$u_0,\ldots,u_b$  (resp. $v_1, \ldots , v_a$)
are the finite (resp. infinite) solutions at~$z=0$ of the equation
$K(z,u)=0$ and the kernel $K$ is defined by~{\rm 
(\ref{kern}).} 
In particular, the generating function for all walks, irrespective of
their endpoint, is
\[
F(z,1)=-\frac{1}{z}\prod_{i=0}^b (1-u_i),
\]
and the \gf \ for {\em excursions}, i.e., walks ending at $0$, is, for $b<0$:
\[
F(z,0)=\frac{(-1)^{b}}{z}\prod_{i=0}^b u_i,
\]
(for $b=0$, the relation becomes $F(z,0)=\frac{(-1)^{b}}{1+z-p_0 z}\prod_{i=0}^b u_i$.)
\end{Proposition}
{These results
 could be derived by a detour via multivariate linear recurrences,
 and the present treatment is closely related to~\cite{BoPe98,Pe98};
 however, our results were obtained independently  in March 1998
 \cite{Ba}.       }

The asymptotic behaviour of the number of $n$-step 
walks can be established via 
singularity analysis or saddle point methods. The series
$u_i$ have ``in general'' a square root singularity, yielding
an asymptotic behaviour  of the form
$A \mu^n n^{-3/2}.$ 
We plan to develop this study in  a forthcoming paper.

\begin{Example}\label{CatalnExample}{\em Catalan numbers}\\
This is the simplest factorial walk,  $(k)\leadsto
 (0)(1)\dots(k)(k+1)$, which corresponds to the ECO-system $(2.f)$.
The  operator $M$ is given by
$$M[f](u)=\frac{f(1)-f(u)}{1-u}+(1+u)f(u).$$
The kernel is $K(z,u)=1-u+z-z(1-u)(1+u)=1-u+zu^2$, hence
$u_0(z)=\frac{1-\sqrt{1-4z}}{2z}$, so that
 $$F(z,1)=-\frac{1-u_0}{z}=\frac{1-2z-\sqrt{1-4z}}{2z^2}=
\sum_{n\ge 1} {2n \choose n}\frac{z^{n-1}}{n+1},$$
the generating function of the Catalan numbers 
(sequence {\bf M1459}\footnote{%
	The numbers {\bf Mxxxx} are identifiers of the
	sequences in {\em The Encyclopedia of Integer Sequences}~\cite{SlPl95}.
}%
).
This result could be expected, given the obvious relation
between these walks and {\L}ukasie\-wicz codes.
\end{Example}

\begin{Example}\label{MotzkinExample}{\em Motzkin numbers}\\
This example, due to Pinzani and his co-authors, is derived from the
previous one by forbidding ``forward'' jumps of length zero. The rule is
then
\[
(k)\leadsto (0)\cdots(k-1)(k+1).
\]
The  operator $M$ is 
$$M[f](u)=\frac{f(1)-f(u)}{1-u}+uf(u).$$
The kernel is $K(z,u)=1-u+z-zu(1-u)=1+z-u(1+z)+zu^2$, leading to
\begin{eqnarray*}
F(z,1)&=&\frac{1-z-\sqrt{1-2z-3z^2}}{2z^2}=1+z+2z^2+4z^3+9z^4+21z^5+O(z^6),
\end{eqnarray*}
the generating function for Motzkin numbers (sequence {\bf
M1184}).
\end{Example}

\begin{Example}\label{SchroederExample}{\em Schr\"oder numbers}\\
This example is also due to the Florentine group.
The rule is $(k) \leadsto (0)\dots (k-1)(k) (k+1)^2$.
{}From Proposition~\ref{thm-algebraic}, we derive
$$F(z,1)=\frac{1-3z-\sqrt{1-6z+z^2}}{4z^2}=1+3z+11z^2+45z^3+197z^4+O(z^5).$$
The coefficients are the Schr\"oder numbers ({\bf M2898}: Schr\"oder's
second problem). We give in Table 1 at the end of the paper a 
 generalization of Catalan and Schr\"oder numbers, corresponding
to the rule $(k) \leadsto (0)\dots (k-1)(k) (k+1)^m$. {This generalized
rule has recently been shown to describe a set of permutations
avoiding certain patterns \cite{Darla}.}
\end{Example}

The above examples were all  quadratic.
However, it is clear  from our treatment that algebraic functions of
arbitrary degree can be obtained: it suffices that the set of
``exceptions'' around $k$ have a span greater than~1.
Let us start with a family of ECO-systems where supplementary forward jumps of
length larger than one are allowed.
\begin{Example}\label{TernaryExample}{\em Ternary trees, dissections of a
polygon, and $m$-ary trees}\\
The ECO-system with axiom $(s_0)=(3)$ and rule \[(k)\leadsto
(3)(4)\cdots(k)(k+1)(k+2)\] is equivalent to the walk \[(k)\leadsto
(0) (1)\cdots(k)(k+1)(k+2).\]
The kernel is $K(z,u)=1-u+zu^3$, and the \gf \ 
$$F(z,1)=\sum_{n\ge1} {3n \choose n} \frac{z^{n-1}}{2n+1}$$
 counts ternary trees  ({\bf M2926}).

More generally, the system with axiom $(m)$
and rewriting rules
\[
(k)\leadsto (m)\cdots(k)(k+1)(k+2)\cdots(k+m-1)
\]
yields the $m$-Catalan numbers, ${mn \choose n}/((m-1)n+1)$, that count $m$-ary
trees. The kernel is $1-u+zu^m$ and the generating function $F(z,1)$
satisfies 
 $F(z,1)=(1+zF(z,1))^m$. 
In particular, the $4$-Catalan numbers ${4n \choose n}/(3n+1)$
appear  in \cite{SlPl95} (sequence {\bf M3587}) and count dissections of a
polygon.

\end{Example}

In the above examples, all backward jumps are allowed. In
other words, each of these examples corresponds to an ECO-system. 
Let us now give an example where backward jumps of length $1$ are
forbidden.

\begin{Example}
{}\\
Consider the following modification of the Motzkin rule:
 \[
(k)\leadsto (0)\cdots(k-2)(k+1).
\]
The kernel is now $K(z,u)=u(1-u)+zu-z(1-u)(u^2-1)$, and, according to
\Ref{biv}, the series $F(z)=F(z,1)$ is given by $F(z)=1/[z(v_1-1)]$, where
$v_1$ satisfies $K(z,v_1)=0$ and is infinite at $z=0$. Denoting
$G=zF(z)$, we find that the algebraic equation defining $G$ is:
$$G=z \ \frac{1+2G+G^2+G^3}{1+G}.$$

\end{Example}

So far, we have only dealt with walks for which the set of allowed
moves was obtained by modifying the interval $[0,k]$ around
$k$. One can also modify this interval around $0$: we shall see -- in
examples -- that the \gf \ remains algebraic. However, it is
interesting to note that in these examples, the kernel method does not
{immediately} provide enough equations between the ``unknown
functions'' to solve the functional equation.

Let us first explain how we modify the interval $[0,k]$ around
$0$. The walks we wish to count are still specified by a multiset $A$
of allowed supplementary jumps and a set $B$ of forbidden backward
jumps. But, in addition, we forbid backward jumps to end up in $C$,
where $C$ is a given finite subset of $\ns$. In other
words, the possible moves from $k$ are given by the rule
$$(k)\leadsto [0, k-1] \setminus(C \cup (k-B)) \ \cup \ (k+A).$$
Again, we can write a functional equation defining $F(z,u)$:
\beq \label{eq-funct-C}
F(z,u)=1+z\left(
\frac{F(z,1)-F(z,u)}{1-u}+P(u)F(z,u)+\sum_{k=0}^{b-1}r_k(u)  F_k(z)
-\sum_{\ga \in C}u^\ga G_\ga (z)\right),
\eeq
where, as above,
$$P(u)=\sum_{\al \in A} u^\al - \sum_{\be \in B} u^{-\be}
\ \ \ \hbox{ and } \ \ \ r_k(u)=\sum_{\be >k, \ \be \in B}u^{k-\be},$$
the new terms in the equations being
$$G_\ga(z)=F(z,1)-\sum_{k=0}^{\ga} F_k(z) - \sum_{\be \in B} F_{\be +
\ga}(z).$$
Observe that the first three terms are the same as in the case
$C=\emptyset$. 
The equation, once multiplied by $u^b(1-u)$, reads $K(z,u)F(z,u)=R(z,u)$
where $K(z,u)$ is given by \Ref{kern} and
$$R(z,u)=u^b(1-u)\left(1+
\frac{zF(z,1)}{1-u}+z\sum_{k=0}^{b-1}r_k(u)  F_k(z)
-z\sum_{\ga \in C}u^\ga G_\ga (z)\right).$$
The kernel is not modified by the introduction of $C$. As above, it
has degree $a+b+1$ in $u$, and admits $b+1$ finite 
roots $u_0, \ldots , u_b$ around $z=0$. However, $R(z,u)$ now involves
$b+1+|C|$ unknown functions, namely $F(z,1)$, the $F_k(z)$, $0\le k <b$
and the $G_\ga(z)$, $\ga \in C$. The degree of $R$ in $u$ is no longer
$b+1$ but $b+c+1$, where $c=\max C$.
The $b+1$ roots of $K$ that can be substituted for $u$ in 
Eq.~\Ref{eq-funct-C} 
provide $b+1$ linear equations between the $b+|C|+1$ unknown
functions. Additional equations will be obtained by extracting the
coefficient of $u^j$ from Eq.~\Ref{eq-funct-C}, for some values of $j$. In
general, we have:
\beq \label{Fj}
F_j(z)=[j=0]+ z\sum_{\al \in A} F_{j-\al}(z)+z[j \not \in C]
\left(F(z,1)-\sum_{k=0}^jF_k(z)-\sum_{\be \in B} F_{j+\be}(z)\right).
\eeq
It is possible to construct a finite subset $S\subset \ns$ such
that the combination of the $b+1$ equations obtained via the kernel
method and the equations \Ref{Fj} written for $j \in S$ determines all unknown
functions as algebraic functions of $z$ -- more precisely, as rational
functions of $z$ and the roots $u_0, \ldots , u_b$ of the
kernel. However, this is a long 
development, and these classes of walks play a marginal role in the
context of ECO-systems. For these reasons, we shall merely give two
examples. {The details
on the general procedure for constructing the set $S$ 
can be found in \cite{mbm}}. 

\begin{Example}
{}\\
This example is obtained by modifying the Motzkin rule of
Example~\ref{MotzkinExample} around {the point} $0$. Take $A=C=\{1\}$ and 
$B=\emptyset$. The rewriting rule is
$$(k) \leadsto (0)(2)(3) \cdots (k-1)(k+1).$$
The functional equation reads
\beq \label{C-ex1}
(1-u+z-zu(1-u))F(z,u)=1-u+zF(z,1)-zu(1-u)G_1(z),
\eeq
with $G_1(z)=F(z,1)-F_0(z)-F_1(z)$. The kernel has a {\em unique\/}
finite root at $z=0$: 
$$u_0=\frac{1+z-\sqrt{1-2z-3z^2}}{2z},$$
{whereas the right-hand side of Eq.~\Ref{C-ex1} contains {\em two\/}
unknown functions.}
Writing Eq.~\Ref{Fj} for $j=0$ and $j=1$ yields
$$F_0(z)=1+z(F(z,1)-F_0(z))\ \ \ \hbox{ and } \ \ \ F_1(z)=zF_0(z).$$
These two equations allow us to express $F_0$ and $F_1$, and hence
$G_1$, in terms of $F(z,1)$:
$$G_1(z)=(1-z)F(z,1)-1.$$
{This equation relates the two unknown functions of Eq.~\Ref{C-ex1}.}
We replace $G_1(z)$ by the above expression in \Ref{C-ex1}, so that
only one unknown function, namely $F(z,1)$, is left. The kernel method
finally gives:
$$F(z,1)=\frac{3-3z^2-2z^3-(1+z)\sqrt{1-2z-3z^2}}{2(1-z-z^2+z^3+z^4)}
=1+z+2z^2+3z^3+6z^4+12z^5+O(z^6).$$
\end{Example}

\begin{Example}
{}\\
Let us choose  $A=\{1\}$, $B=\{2\}$ et $C=\{2\}$.
The rewriting rule is now:
$$(k) \rightarrow (0)(1)(3)(4)(5) \ldots (k-3)(k-1)(k+1).$$
The functional equation reads
$$\left[u^2(1-u)+zu^2-zu^3(1-u)+z(1-u)\right]F(z,u)\hskip 8cm $$
\beq \hskip 2cm =u^2(1-u)+
zu^2F(z,1)+z(1-u)\left[F_0(z)+uF_1(z)\right] -zu^4(1-u) G_2(z)
,\label{eq2}\eeq
with $G_2(z)=F(z,1)-F_0(z)-F_1(z)-F_2(z)-F_4(z)$. 
Only three roots, $u_0, u_1, u_2$ can be substituted for $u$ in the
kernel, while the right-hand side of the equation contains four
unknown functions, $F(z,1), F_0(z), F_1(z)$ and $G_2(z)$.
Writing \Ref{Fj} for $j=0,1$ and $2$ yields
$$\begin{array}{ll}
F_0(z)& =1+z\left[ F(z,1)-F_0(z)-F_2(z)\right], \\
F_1(z) &= zF_0(z) + z\left[ F(z,1)-F_0(z)-F_1(z)-F_3(z)\right],\\
F_2(z)&=zF_1(z).
\end{array}$$
The second equation is not of much use but, by combining the first and
third one, we find
$$F_0(z) =\frac{1+z\left[ F(z,1)-zF_1(z)\right]}{1+z}.$$
Replacing  $F_0(z)$ by this expression in  \Ref{eq2} gives:
$$\left[u^2(1-u)+zu^2-zu^3(1-u)+z(1-u)\right]F(z,u)
 =u^2(1-u)+\frac{z(1-u)}{1+z}\hskip 4cm $$
\beq +zF(z,1) \left[
u^2+\frac{z(1-u)}{1+z}\right] +z(1-u)F_1(z) \left[
u-\frac{z^2}{1+z}\right] -zu^4(1-u)G_2(z) .\label{eq2bis}\eeq
We are left with three unknown functions, related by three linear
equations obtained by cancelling the kernel. Solving these equations
would give $F(z,1)$ as an enormous rational function of $z$, $u_0,u_1$
and $u_2$, symmetric in the $u_i$. This implies that $F(z,1)$ can also
be written as a rational function of $z$ and $v\equiv v_1$, the
{fourth and last}
root of the kernel. In particular, $F(z,1)$ is algebraic of degree at
most $4$. 

In order to obtain directly an expression  of $F(z,1)$ in terms of $z$
and $v$, we can proceed as follows. Let $R'(z,u)$ denote the
right-hand side of Eq.~\Ref{eq2bis}. Then $R'(z,u)$ is a polynomial in
$u$ of degree $5$, and three of its roots are $u_0, u_1,
u_2$. Consequently, as a polynomial in $u$,
the kernel $K(z,u)$ divides $(u-v)R'(z,u)$. 

Let us evaluate $(u-v)R'(z,u)$ modulo $K(z,u)$: we obtain a
polynomial of degree $3$ in $u$, whose coefficients depend on $z, v,
F(z,1), F_1(z)$ and $G_2(z)$. This polynomial has to be zero: this
gives a system of four (dependent) equations relating the three unknown
functions $F(z,1), F_1(z)$ and $G_2(z)$. Solving the first three of
these equations yields
\begin{eqnarray*}
F(z,1)&=&\frac {1+z+z^2-(z+1)zv+(z+1)zv^2-z^2v^3}
{1-z^2-z(1-z^2)v+ z^3v^3}\\
&=& 1+z+2z^2+3z^3+6z^4+11z^5+23z^6+47z^7+101x^8+O(z^9).\\
\end{eqnarray*}
 Eliminating $v$ between this expression and $K(z,v)=0$ gives
a quartic equation satisfied by $F(z,1)$.
\end{Example}
\section{Transcendental systems}
\label{SectionTranscendant}

\subsection{Transcendence}
The radius of convergence of an algebraic series is always positive.
Hence,  one possible reason for a system to give a transcendental series is the
fact that its coefficients grow too fast, so that its radius of convergence is
zero. This is the case for 
System~$(2.h)$, as proved by the following proposition.

\begin{Proposition} \label{zero1}
Let $b$ be a nonnegative integer. 
For  $k \ge 1$, let $m(k)=|\{i: \ e_i(k) \ge k-b\}|$.
Assume that:

$1.$  for all $k$, there exists a forward
 jump from~$k$ (i.e., $e_i(k)>k$ for some $i$), 

$2.$ the sequence  $(m(k))_k$ is  nondecreasing and tends to infinity.

\noindent  Then the {(ordinary)} generating function of the system has radius of
convergence $0$.
\end{Proposition}
\Proof
Let $s_0$ be the axiom of the system.
{Let us denote by $h_n$ the product} $m(s_0+b)m(s_0+2b)\cdots m(s_0+nb)$. Let us
prove that the generating tree contains  at least  $h_n$ nodes at
level $n(b+1)$. 
At level $nb$, take a node $v$ labeled $k$, with $k \ge
s_0+nb$.
Such a node exists thanks to the first assumption.
 By definition of $m(k)$, this
node $v$ has  $m(k)$ sons whose label is at least $k-b$. As $m$
is non decreasing, $v$ has at least $m(s_0+nb)$ sons of label at least
$s_0+(n-1)b$. Iterating this procedure shows that, at level
$nb+i$, at least $m(s_0+(n-i+1)b)\cdots m(s_0+nb)$ descendants
of $v$ have a label larger than or equal to $s_0+(n-i)b$, for $0<
i \le n $.  In particular, for $i=n$, we obtain at level $n(b+1)$ at
least $h_n$ descendants of $v$ whose label is at least $s_0$.

Hence  $f_{n(b+1)} \ge h_n$. But as $h_{n}/h_{n-1} = m(s_0+nb)$
goes to infinity with $n$, the series $\sum_n h_n z^{n(b+1)}$ has radius of
convergence $0$, and the same is true for $F(z)=\sum_n f_nz^n$.
\qed

In particular, this proposition implies that {\em the \gf \ of any
ECO-system in which the length of backward jumps is bounded has
radius of convergence $0$.} Many  examples of this type will be given
in the next subsection, in which we shall study whether the
corresponding \gf \ is holonomic or not.
  The following example, in which backward jumps are not bounded,
 was suggested by Nantel Bergeron.

\begin{Example}
{\em A fake factorial walk} \\
Consider the system with axiom $(1)$ and rewriting rules
$\{(k) \leadsto (2)(4) \cdots (2k)\}$.  Proposition~\ref{zero1}
applies with $b=0$ and $m(k)=1+ \lfloor k /2 \rfloor$. Note that 
the radius of convergence of $F(z)$ is zero 
although {\em all\/} the functions $e_i$ are bounded, and indeed
constant: $e_i(k)=2i$ for all $k \ge i$.
The series $F(z)$ is of course 
transcendental. Note, however,
that $F(z,u)$ satisfies a functional equation that is at first sight
reminiscent of the equations studied in Section \ref{SectionAlgebraic}:
$$F(z,u)=u+ zu^2 \ \frac{F(z,1)-F(z,u^2)}{1-u^2}.$$
\end{Example}

The following example shows that Proposition~\ref{zero1} is not far
from optimal: an ECO-system in which all 
functions $e_i$ grow linearly can have a finite radius of convergence.

\begin{Example} 
{}\\
The system with axiom $(1)$ and rules $(k) \leadsto (\lceil k/2
\rceil)^{k-1}(k+1)$ leads to a \gf \ with a positive radius of
convergence.\\
Let us start from the recursion defining the numbers $f_{n,k}$. We
have $f_{0,1}=1$ and for $n \ge 1$,
$$f_{n+1,k}=f_{n,k-1}+(2k-1)f_{n,2k}+(2k-2)f_{n,2k-1}.$$
 The largest label occurring at level $n$ in the tree is $n+1$. Let us
introduce the numbers $g_{n,k}=f_{n,n-k+1}$, for $k \le n$. The above
recursion can be rewritten as:
\beq g_{n+1,k}=g_{n,k}+(2n-2k+3)g_{n,2k-n-3}+(2n-2k+2)g_{n,2k-n-2}.
\label{g} \eeq
We have $g_{n,k}=0$ for $k<0$. Hence Eq.~\Ref{g} implies that for
$k \ge 0$, the sequence $(g_{n,k})_{n}$ is nondecreasing and reaches a
constant value $g(k)$ as soon as $n\ge 2k-1$ {(see Table~\ref{fgtable}).}

\begin{table}
$$\begin{array}{l|ccccccc}
	n \ k& 1 	&2	&3	&4	&5	&6		\\
\hline
0	&\bf{1} \\
1	&\bf{0}	&\bf{1}	\\
2	& 1	&\bf{0}	&\bf{1}	\\
3 	&0	&\bf{3}	&\bf{0}	&\bf{1}	\\
4 	&3	&3	&\bf{3}	&\bf{0}	&\bf{1}	\\
5 	&3 	&9	&\bf{7}	&\bf{3}	&\bf{0}	&\bf{1}	\\
\end{array}
\hskip 2cm
\begin{array}{l|ccccccc}
	n \ k& 0	&1 	&2	&3	&4	&5	\\
\hline
0	&\bf{1} \\
1	&\bf{1}	&\bf{0}	&\\
2	&\bf{1}	&\bf{0}	& 1	\\
3 	&\bf{1}	&\bf{0}	&\bf{3}	&0\\
4 	&\bf{1}	&\bf{0}	&\bf{3}	&3	&3	\\
5 	&\bf{1}	&\bf{0}	&\bf{3}	&\bf{7}&9	&3	\\
\end{array}
$$
\caption{The numbers $f_{n,k}$ and  $g_{n,k}$. Observe the convergence
of the coefficients.}
\label{fgtable}
\end{table}
Going back to the number $f_n$ of nodes at level $n$, we have $$f_n
=\sum_{k=0}^n g_{n,k} \le \sum_{k=0}^n g(k).$$
But $$\sum_{n\ge 0}z^n  \sum_{k=0}^n g(k) = \frac{1}{1-z}
\sum_{k=0}^n g(k)z^k,$$
and hence it suffices to prove that the \gf \ for the numbers $g(k)$
has a finite radius of convergence, that is, that these numbers grow
at most exponentially.

Writing \Ref{g} for $n+1=2k-i$, for $1\le i\le k$, we obtain:
$$g_{2k-i,k}=g_{2k-i-1,k}
+(2k-2i+1)g_{2k-i-1,i-2}+(2k-2i)g_{2k-i-1,i-1}.$$
Iterating this formula for $i$ between $1$ and $k$ yields
\begin{eqnarray*}
g(k) &=& g_{2k-1,k}= \sum_{i=1}^k \left[
(2k-2i+1)g_{2k-i-1,i-2}+(2k-2i)g_{2k-i-1,i-1}\right]\\
&\le & \sum_{i=1}^k \left[
(2k-2i+1)g(i-2)+(2k-2i)g(i-1)\right]=
 \sum_{i=0}^{k-2} (4k-4i-5)g(i).
\end{eqnarray*}
This inequality, together with the fact that $g(0)=1$, implies that
for all $k \ge 0$, $g(k)\le \tilde g(k)$, where the sequence $\tilde
g(k)$ is defined by $\tilde g(0)=1$ and $\tilde g(k) =
\sum_{i=0}^{k-2} (4k-4i-5)\tilde g(i)$ for $k >0$. But the series
$\sum_k \tilde g(k) z^k$ is rational, equal to
$(1-z)^2/(1-2z-2z^2-z^3)$, and has a finite radius of
convergence. Consequently, the numbers $\tilde g(k)$ and $ g(k) $
grow at most exponentially.
\end{Example}

\medskip
Algebraic generating functions are strongly constrained in their
algebraic structure (by a polynomial equation) as well as in their
analytic structure (in terms of singularities and asymptotic
behaviour). In particular, they have a finite number of
 singularities, which are algebraic numbers, and they admit local asymptotic
expansions that  involve only rational exponents.
{\em A contrario\/}, a generating function that has infinitely many
singularities ({\em e.g.}, a natural boundary) or that involves a
transcendental element ({\em e.g.}, a logarithm) in a local asymptotic
expansion is by necessity transcendental; see~\cite{Flajolet87}
for a discussion of such transcendence criteria. In the case of
generating trees, 
this means that the presence of a condition involving a transcendental
element is expected to lead to a transcendental generating function.
This is the case in the following example.

\begin{Example}{\em A Fredholm system}\label{FredholmExample}\\
 We examine System~$(2.e)$, in which
the rules are irregular at powers of~2:
$$ (s_0)=(2), \ \ \  \ (k) \leadsto (2)^{k-2} (3 - [\exists p\!:\!k=2^p])
(k+1), \ \ \ k\ge 2.$$
This example will involve the Fredholm series
$h(z):=\sum_{p\ge 1} z^{2^p}$, which is well-known to admit the unit
circle as a natural boundary. (This can be seen by way of the
functional equation $h(z)=z^2+h(z^2)$, from which there results that
$h(z)$ is infinite at all iterated square-roots of unity.)
According to Eq.~\Ref{equationF_k}, we have, for $k>3$,
$F_k(z)=zF_{k-1}(z),$
so that $$F_k(z)=z^{k-3}F_3(z) \ \ \ \hbox{ for } k\ge 3.$$
Now, writing Eq.~\Ref{equationF_k} for $k=2$ gives
\begin{eqnarray*}
F_2(z)&=&1+z\sum_{k\ge 3} (k-2)F_k(z) +z \sum_{p\ge 1} F_{2^p}(z)\\
&=& 1+\frac{z}{(1-z)^2} F_3(z) +z F_2(z) + F_3(z) \left(
\frac{h(z)}{z^2}-1\right)\\
&=& 1+zF_2(z) + F_3(z) \left( \frac{z}{(1-z)^2} +
\frac{h(z)}{z^2}-1\right).
\end{eqnarray*}
For $k=3$, we obtain:
\begin{eqnarray*}
F_3(z)&=&zF_2(z)+z\sum_{k\ge 3, \ k \not = 2^p} F_k(z)\\
&=& zF_2(z)+F_3(z)\left( \frac{1}{1-z}-\frac{h(z)}{z^2}\right).
\end{eqnarray*}
Solving for $F_2(z)$ and $F_3(z)$, then summing ($F(z) = F_2(z) +
F_3(z) / (1-z)$), we obtain:
\[
F(z) = \frac {  (1-z)^2 h(z)} 
{ (1-2z) (1-z)^2 h(z) - z^4  }=1+2z+5z^2+14z^3+39z^4+108z^5+O(x^6)
.
\]
The functions $h(z)$ and $F(z)$ are rationally related, so that 
$F(z)$ is itself transcendental. The series $h$ has radius $1$, but
the denominator of $F$ vanishes before $z$ reaches $1$ -- actually,
before $z$ reaches $1/2$. Hence the radius of $F$ is the smallest root
of its denominator. Its value is easily
determined numerically and found to be about 0.360102.
\end{Example}

\subsection{Holonomy}
In the transcendental case, one can also discuss the {\em
holonomic\/} character of the generating function $F(z)$.

A series is said to be 
{\em holonomic\/}, or {\em D-finite\/} \cite{stanley}, if it satisfies
a linear differential equation with polynomial coefficients in
$z$. Equivalently, its coefficients $f_n$ satisfy a linear recurrence
relation with polynomial coefficients in $n$. Consequently, given a
sequence $f_n$, the ordinary \gf \ $\sum_n f_n z^n$ is holonomic if
and only if the exponential \gf \ $\sum_n f_n z^n /n!$ is holonomic.
The set of holonomic series has nice closure properties: the sum or product of
two of them is still holonomic, and the substitution of an algebraic
series into an holonomic one gives an holonomic series.
  Holonomic series include algebraic series, and have a
finite number of singularities. 
This implies that Example~\ref{FredholmExample}, for which $F(z)$ has
a natural boundary, is not holonomic.

We study below five ECO-systems that, at first sight, do not
look to be very
different. In particular, for each of them, forward
and backward jumps are bounded. Consequently, Proposition~\ref{zero1}
implies that the corresponding ordinary \gf \ has radius of convergence
zero. However, we shall see that the first three systems have an
holonomic \gf , while the last two have not. We have no general
criterion that would allow us to distinguish between systems leading
to holonomic generating functions and those leading to 
nonholonomic generating functions.

Among the systems with bounded jumps, those for which $e_i(k)-k$
belongs
to $\{-1,0,1\}$ for all $i\le k$ 
have a nice property: the \gf \ for the corresponding {\em
excursions\/} (walks starting and ending at level $0$) can be written
as the following continued fraction \cite{flajolet-fracont}:
$$\frac{1}{\displaystyle 1-b_0z-
\frac{\displaystyle a_1c_0z^2}{\displaystyle 1-b_1z-
\frac{\displaystyle a_2c_1z^2}{\displaystyle 1-b_2z-
\frac{\displaystyle a_3c_2z^2}{\cdots}}}},$$
where the coefficients $a_k, b_k$ and $c_k$ are the multiplicities appearing
in the rules, which read $(k) \leadsto (k-1)^{a_k}(k)^{b_k}(k+1)^{c_k}.$

\begin{Example} \label{Arrangements}{\em Arrangements}\\
The system $(k) \leadsto (k)(k+1)^{k-1}$ with axiom~$(s_0)=(2)$ generates a
sequence that starts with $1,2,5,16,65, 326$ 
({\bf M1497}). It is not hard to see that the 
triangular array $f_{n,k+2}$ is given by the arrangement
numbers $k!{n \choose k}$, so that the {\em exponential\/} generating
function (EGF) of 
the sequence is $$\tilde F(z,u) 
= \sum_{n\ge 0,k\ge 2} f_{n,k}u^k \frac{z^n}{n!}
=\frac{u^2 e^z}{1-uz}.$$ This system satisfies the
conditions of Proposition~\ref{zero1} with $b=0$ and
$m(k)=k$. Accordingly, one has $f_n\sim e\, n!$, so that the {\em ordinary\/}
generating function $F(z)$  has radius of convergence~0 and  cannot be
algebraic. However, $\tilde F(z,1)=e^z/(1-z)$ is holonomic, and so is
$F(z)$. 
\end{Example}

 \begin{Example} \label{ExampleInvolutions}{\em Involutions and
Hermite polynomials}\\
The system $(k) \leadsto (k-1)^{k-1}(k+1)$ with axiom~$(s_0)=(1)$
 generates a sequence that starts with $1,1,2,4, 10,26, 76$ ({\bf
 M1221}). These numbers count involutions: more precisely, one easily
derives  from the recursion satisfied by the coefficients $f_{n,k}$ that
$f_{n,k}$ is the number of involutions on $n$ points, $k-1$ of which
are fixed.  Proposition~\ref{zero1} applies with $b=1$ and $m(k)=k$.

 The corresponding EGF is
\beq \tilde F(z,u)= \sum_{n\ge 0,k\ge 1} f_{n,k}u^k \frac{z^n}{n!}
= u\exp\left (zu+\frac{z^2}{2}\right),\label{inv-egf} \eeq
and its value at $u=1$ is holonomic.

The polynomials $f_n(u)=\sum_k f_{n,k}u^k$ counting involutions on $n$
points are in fact closely related to the Hermite polynomials, defined by:
$$\sum_{n\geq 0} H_n(x) \frac{t^n}{n!}=\exp\left(xt-\frac{t^2}{2}\right).$$
Indeed, comparing the above identity with \Ref{inv-egf} shows that
$ f_n(u)=u\ i^nH_n(-i u).$
\end{Example}

\begin{Example}{\em Partial permutations and Laguerre polynomials}
\label{ExampleLaguerre}\\
The rewriting rule $(k)\leadsto  (k+1)^{k-1}
(k+2)$, taken with axiom $(2)$, generates a sequence that starts with
$1,2,7,34,209,...$ ({\bf M1795}). From the recursion satisfied by the
coefficients $f_{n,k}$, we derive that $f_{n,n+k}$  is the number of
partial injections of $\{1,2,\ldots , n\}$ into itself in which $k-2$
points are unmatched.
From this, we obtain:
$$\tilde F(z,u)=
\frac{u^2}{1-uz}\exp\left( \frac{u^2z}{1-uz}\right)=u^2\sum_{n\geq 0}
L_n(-u)  \frac{(uz)^n}{n!} $$
where $L_n(u)$ is the $n$th Laguerre polynomial.
Again, $\tilde F(z,1)$ is holonomic.
\end{Example}

\noindent The next two systems, as announced,
lead to nonholonomic generating functions.

\begin{Example}\label{ExampleStirling}{\em Set partitions and 
 Stirling polynomials}\\
Let us consider the system $[(1),(k)\leadsto (k)^{k-1}(k+1)]$. From
the recursion satisfied by the coefficients $f_{n,k}$, we derive that
$f_{n,k+1}$ is equal to the Stirling number of
 the second kind  $\stir{n}{k}$, which counts partitions of  $n$ objects
into $k$ nonempty subsets. The corresponding EGF is
$$\tilde F(z,u) = 
u \exp\left(u(\exp{z}-1)\right).$$ 
At $u=1$, this \gf  \ specializes to 
$$\tilde
 F(z,1)=\exp(\exp(z)-1))=\sum_{n \ge 0} B_n \frac{z^n}{n!}=
1+z+2\frac{z^2}{2!}+5
\frac{z^3}{3!}+15\frac{z^4}{4!}+52\frac{z^5}{5!}+203\frac{z^6}{6!}+\dots$$
This is the exponential generating function of the Bell
 numbers ({\bf M1484}).  It is known that $\log B_n =n \log n -n \log
\log n +O(n)$ (see~\cite{Odlyzko}), and this 
cannot be the asymptotic behaviour of the logarithm of the
coefficients of an holonomic 
series  (see~\cite{WiZe85b} for admissible types).
Hence, $\tilde
F(z,1)$, as well as $F(z,1)$, is nonholonomic.
\end{Example}

\begin{Example}\label{ExampleBessel}{\em Bessel numbers}\\
We study the system with axiom~(2) and rewriting rules
\beq 
(2)\leadsto(2)(3),\ \ \  (k)\leadsto (k-1)(k)^{k-2}(k+1), \ \  k\ge
3. 
\label{bess}
\eeq
We shift the labels by $2$ to obtain  a walk model with
axiom $(0)$ and rules
$$ (0)\leadsto (0)(1), \ \ \  (k)\leadsto (k-1)(k)^k(k+1), \ \ k\ge
1.$$
 The corresponding bivariate \gf \ $F(z,u)$ 
satisfies the functional differential equation
\[
F(z,u)\left(1-z(u+u^{-1})\right)
=1+z(1-u^{-1})F(z,0)+zu\frac{\partial F}{\partial u}(z,u),
\]
which is certainly not obvious to solve.
However, as observed in  \cite{flajolet-fracont}, it is easy to obtain
a continued fraction expansion of the excursion \gf :
\[
F(z,0)=1+z+2z^2+4z^3+9z^4+\cdots=\cfrac{1}{1-z-
	\cfrac{z^2}{1-z-
	 \cfrac{z^2}{1-2z-
	  \cfrac{z^2}{1-3z-
		\ddots
	}}}}=\frac{1}{1-z-z^2B(z)},
\]
where $B(z)=\sum_n B^*_nz^n=1+z+2z^2+5z^3+14z^4+43z^5+143z^6+\cdots
$
is the generating function of Bessel numbers ({\bf M1462}) and counts
non-overlapping partitions \cite{FlSc90}. 
As $F(z,0)$ itself, the series $B(z)$ has radius of convergence
zero. The fast increase of $B^*_n$ entails 
\[
[z^n]F(z,0)\sim B_{n-2}^*.
\] {}From~\cite{FlSc90}, we know that $\log B^*_n = n \log n -n \log
\log n +O(n).$
Again,  this prevents $F(z,0)$ from being holonomic.

In order to prove that $F(z,1)$ itself is nonholonomic, we are going
to prove that its coefficients $f_n$ have the same asymptotic
behaviour as the coefficients of $F(z,0)$. Clearly,
$$[z^n]F(z,0) =f_{n,0} \le\sum_k f_{n,k} =f_n.$$
To find an upper bound for $f_n$, we compare the system~\Ref{bess}
(denoted $\Sigma_1$ below) to
the system $\Sigma_2$ with axiom $(2)$ and rule  $(k)\leadsto
(k)^{k-1}(k+1)$. This system generates a tree with counting sequence
$g_n$. The form of the rules implies that the (unlabeled) tree associated
with $\Sigma_1$ is a subtree of the tree associated with $\Sigma_2$. Hence
$f_n \le 
g_n$. Comparing $\Sigma_2$ to the system studied in the previous example
shows that $g_n$ is the Bell number $ B_{n+1}$, the logarithm of which
is also known to be  $n \log n -n \log
\log n +O(n)$ (see~\cite{Odlyzko}). Hence  $\log f_n = n \log n -n \log
\log n +O(n),$ and this prevents the series  $F(z,1)$ from being holonomic.
\end{Example}

\newpage
\noindent {\bf{A small catalog of ECO-systems}}\\
To conclude, we present in Table~2 a small catalog of ECO-systems that
lead to sequences of 
combinatorial interest. Several examples are detailed in the paper;
others are due to West \cite{West95,West96} or Barcucci, Del Lungo,
Pergola, Pinzani \cite{ECOArticle,BaDePePiFPSAC98,BaDePePi,BaDePePi_2}, 
or are folklore.  Each of them is an instance of
application of our criteria.

\begin{table}[t]
\begin{small}
\begin{tabular}{|c|l|l|c|c|}
\hline
Axiom & System & Name & Id. & Generating Function \\
\hline
&\bf{Rational OGF}  &&& {\bf OGF}\\
$(1)$ & $(k) \leadsto (k)^{k-1} ((k \bmod 2)+1)$ &
Ex. \ref{FibonacciExample}: Fibonacci & M0692 & $\frac{1}{1-z-z^2}$ 
\makebox(0,12)[0cm]{}
\\
$(2)$ & $(k) \leadsto (2)^{k-1} (k+1)$ & Ex. \ref{FiboBisectionExample}: even Fibonacci & M1439 & $\frac{1-z}{1-3z+z^2}$ 
\makebox(0,12)[0cm]{}
\\
$(3)$ & $(k) \leadsto (2)^{k-1} (k+1)$ &
Ex. \ref{FiboBisectionExample}: odd Fibonacci & M2741 & $\frac{1}{1-3z+z^2}$ 
\makebox(0,12)[0cm]{}
\\
&&&&\\
\hline
&\bf{Algebraic OGF} &&&  {\bf OGF}\\
$(1)$ & $(k) \leadsto (1) 
\cdots (k-1) (k+1)$ & Ex. \ref{MotzkinExample}: Motzkin numbers & M1184 &  
$\frac{1-z-\sqrt{1-2z-3z^2}}{2z^2}$
\makebox(0,12)[0cm]{}
\\
$(2)$ & $(k) \leadsto (2) 
\cdots (k) (k+1)$ & Ex. \ref{CatalnExample}: Catalan numbers & M1459 &  
$\frac{1-2z-\sqrt{1-4z}}{2z^2}$
\makebox(0,12)[0cm]{}
\\
$(3)$ & $(k) \leadsto (3) 
\cdots (k) (k+1)^2$ & Ex. \ref{SchroederExample}: Schr\"oder numbers & M2898 &  
$\frac{1-3z-\sqrt{1-6z+z^2}}{4z^2}$
\makebox(0,12)[0cm]{}
\\
$(4)$ & $(k) \leadsto (4) 
\cdots (k) (k+1)^3$ &\hskip 20mm --- & M3556 &  
$\frac{1-4z-\sqrt{1-8z+4z^2}}{6z^2}$
\makebox(0,12)[0cm]{}
\\
$(m)$ & $(k) \leadsto (m) 
\cdots (k) (k+1)^{m-1}$ &\hskip 20mm --- &--- &  
$\frac{1-mz-\sqrt{1-2mz+(m-2)^2z^2}}{2(m-1)z^2}$
\makebox(0,12)[0cm]{}
\\
$(3)$ & $(k)\leadsto (3) 
\cdots 
(k+2)$ & Ex. \ref{TernaryExample}: Ternary trees & M2926 &  
$F=(1+zF)^3$ 
\makebox(0,12)[0cm]{}
\\
$(4)$ & $(k)\leadsto (4)
\cdots 
(k+3)$ & Ex. \ref{TernaryExample}: Dissections of a polygon & M3587 & 
$F=(1+zF)^4$ 
\\
$(m)$ & $(k)\leadsto (m)\cdots
(k+m-1)$ & Ex. \ref{TernaryExample}: $m$-ary trees & &  
$F=(1+zF)^m$ 
\\
&&&&\\
\hline
&\bf{Holonomic}  &&& {\bf EGF}\\
&{\bf  transcendental OGF} &&&\\

$(1)$ & $(k) \leadsto (k+1)^k$ & Permutations & M1675 &  $1/(1-z)$ \\

$(2)$ & $(k) \leadsto (k)(k+1)^{k-1}$   & Ex. \ref{Arrangements}: 
Arrangements & M1497 & $e^z/(1-z)$ \\

$(1)$ & $(k) \leadsto (k-1)^{k-1} (k+1)$  &
Ex. \ref{ExampleInvolutions}: Involutions & M1221 &
$e^{z+\frac{1}{2}z^2}$\\

$(2)$ & $(k) \leadsto (k+1)^{k-1} (k+2)$ & Ex.~\ref{ExampleLaguerre}:
Partial permutations &
M1795 & $e^{z/(1-z)}/(1-z)$ \\

$(2)$ & $(k) \leadsto (k-1)^{k-2} (k) (k+1)$  & Switchboard problem & M1461 & $e^{2z+\frac{1}{2}z^2}$\\

$(2)$ & $(k) \leadsto (k-1)^{k-2} (k+1)^2$  & Bicolored involutions  & M1648 & $e^{2z+z^2}$\\

&&&&\\
\hline
&\bf{Nonholonomic OGF} && & {\bf EGF}\\
$(1)$ & $(k) \leadsto (k)^{k-1} (k+1)$  & Ex. \ref{ExampleStirling}: 
Bell numbers & M1484 &
$e^{e^z-1}$\\
$(2)$ & $(k) \leadsto (k)^{k-2} (k+1)^2$ & Bicolored partitions  & M1662
& $e^{2(e^z-1)}$\\
$(2)$ & $(k)\leadsto (k-1)(k)^{k-2}(k+1)$ & Ex. \ref{ExampleBessel}: Bessel numbers & M1462
& --- \\
\hline
\end{tabular}
\end{small}
 \caption{ Some ECO-systems of combinatorial interest.}
\end{table}

\bigskip\noindent
{\bf Acknowledgements.}
We thank Elisa Pergola and Renzo Pinzani who presented us the problem
we deal with in this paper. We are also very grateful for helpful
discussions with Jean-Paul Allouche.

\newpage

\bibliographystyle{plain}

\bibliography{eco}

\newpage
\bigskip

\noindent {\em Cyril Banderier, Philippe Flajolet} \\
Projet Algorithmes\\
 INRIA Rocquencourt\\
F-78153 Le Chesnay\\
 {\sc France}\\
{\tt Cyril.Banderier@inria.fr}\\
{\tt Philippe.Flajolet@inria.fr}
\bigskip 

\noindent {\em Mireille Bousquet-M\'elou}\\
LaBRI, Universit\'e Bordeaux 1\\
351 cours de la Lib\'eration\\
F-33405 Talence Cedex\\
 {\sc France}\\
{\tt bousquet@labri.u-bordeaux.fr}

\bigskip 

\noindent {\em Alain Denise, Dominique Gouyou-Beauchamps}\\
LRI, B\^atiment 490\\
Universit\'e Paris-Sud XI \\
F-91405 Orsay Cedex\\
 {\sc France}\\
{\tt Alain.Denise@lri.fr}\\
{\tt dgb@lri.fr}

\bigskip 

\noindent {\em Dani\`ele Gardy}\\
Universit\'e de Versailles/Saint-Quentin \\
Laboratoire PRISM\\
45, avenue des \'Etats-Unis\\
F-78035 Versailles Cedex\\
 {\sc France}\\
{\tt Daniele.Gardy@prism.uvsq.fr}

\end{document}